\documentclass[11pt,centertags]{amsart}
\usepackage{amsmath, amssymb, amscd}
\font\bbb=msbm10

\newtheorem{theorem}{Theorem}[section]
\newtheorem{proposition}{Proposition}[section]
\newtheorem{lemma}{Lemma}[section]

\def\ssu{\subset}

\def\sm{\setminus}

\def\zz{ \hbox{{\bbb Z}}}

\def\nn{ \hbox{{\bbb N}}}
\def\rr{ \hbox{{\bbb R}}}

\def\De{ \Delta}

\def\om{ \omega}
\def\Om{ \Omega}

\def\si{ \sigma}

\def\ep{ \epsilon}
\def\til{ \tilde}
\def\<{\langle}
\def\>{\rangle}

\begin{document}

\title{Amenability of Universal $2$-Grigorchuk group}
\author[R. Muchnik]{R. Muchnik}\thanks{The author was supported by NSF Postdoctoral Fellowship grant DMS-0202457}
\address{University of Chicago} 
\email{roma@math.uchicago.edu} 

\begin{abstract} We consider the universal Grigorchuk $2$-group, i.e., the group such that every Grigorchuk $2$-group is a quotient.
We show that this group has a nice universal representation in the
group of all functions Func($\{0,1,2\}^{\nn}, Aut(T_2))$, where $T_2$
is a group of automorphism of the binary tree. Finally, we prove
that this universal Grigorchuk $2$-group is amenable.  The proof is
an application of the ``M\"{u}nchhausen trick'' developed by V.
Kaimanovich.
\end{abstract}
\maketitle

\section{Introduction.}
Consider a class of Grigorchuk $2$-groups $G_{\om}$.  These groups
are parametrized by infinite sequences $\om \in \{0,1,2\}^{\nn}$.
They have subexponential growth. They have a nice generating set,
i.e., they are all generated by $4$-elements $a_{\om},
b_{\om},c_{\om},d_{\om}$. In this short note I would like to
consider a universal Grigorchuk $2$-group $Gr_2$, which is a
quotient of the free group $F_4=\<a,b,c,d\>$ by a set of words in
$a,b,c,d$ such that these word are trivial for every Grigorchuk
group $G_{\om}$.

Let $\pi_{\om}:F_4\to G_{\om}$ be a map sending generator $a, b,c,d$ to $a_{\om},b_{\om},c_{\om},d_{\om}$.  Then  $Gr_2$ is defined as
$$Gr_2= F_4/ (\cap_{\om} Ker(\pi_{\om}))$$

I will show that $Gr_2$, has a nice representation in the group of
all functions from  the set $\{0,1,2\}^{\nn}$ to the group of
automorphisms of the binary tree $Aut(T_2)$.  This group is
self-similar.  It is self-similar with respect to an isomorphic
group. Because of this technicality I will introduce the notion of
a nested structure.  This is essentially an embedding of $Gr_2$ to
the matrix group on the group algebra of $Gr_2$.

The main result of this note is
\begin{theorem} \label{thm:main} $Gr_2$ is amenable \end{theorem}

To prove this I  invoke a method suggested by V. Kaimanovich
\cite{K04} called ``M\"{u}nchhausen trick''.  The proof is quite simple and follows along the line of \cite{K04}.

\section{Aknowledgements.}
I would like to thank Anna Erschler and Volodya Nekrashevich for
suggesting the problem. Also I am especially grateful to Vadim
Kaimanovich for his help and his preprint, parts of which are
reproduced in this note. The majority of the work on this note was
done during my stay in Bremen at International University of Bremen,
which I thank for their warm hospitality.

\section{Nested structures.}
Let $H$ be a group and $X$ be a set. Then the permutation wreath
product of $H$ and $X$ is $H^X \leftthreetimes Sym(X)$, where the
product is defined as
$$(g^{x}, \si)(h^{x}, \tau) = ((gh^{\si})^{x}, \si\tau).$$

{\bf Definition:} We say $G$ is nested inside $H$ with finite set
$X$ if there exists a group homomorphisms embedding $\phi: G \to
{\bf\til{G}}=H^{X}\leftthreetimes Sym(X),$ such that the image of
the induced map $P_{\phi}:G\stackrel{\phi}{\rightarrow}
{\bf\til{G}}\to Sym(X)$ is irreducible. We denote this structure
$(G,H,X,\phi)$.

\medskip

Now if $G$ is nested in $H$ with set $X$,  I abuse the notation and
identify $G$ with its image and write $g=({\bf h}_g, \si_g) \in
{\bf\til{G}}$, where ${\bf h}_g \in H^X$ and $\si_g \in Sym(X)$.

\section{Random walks on groups with nested structures.}
Next two sections are taken directly from \cite{K04}. However I
reproduce them here for the reader's convenience, and because I
will need to make only a few tiny modifications.

\subsection{Preliminaries.}
Recall that the (right) random walk on a countable group $G$
determined by a probability measure $\mu$ is the Markov chain with
the state space $G$ and the transition probabilities $$\pi_g(k)=
\mu(g^{-1}k)$$ which are equivariant with respect to the left action
of the group on itself. In other words, for a point $g$ the random
walk moves at the next moment of time to the point $gk$, where the
random increment from $k$ is chosen according to the distribution
$\mu$. We shall use for this description of transition probabilities
of the random walk $(G, \mu)$ the notation
\begin{align} \label{form:notation}g_{\stackrel{\longrightarrow}{k \sim \mu}} gk. \end{align}

Thus, if the random walk starts at the moment $0$ from a point
$g_0$, then its position at time $n$ is
$$g_n=g_0k_1k_2\dots k_n,$$
where $k_i\in G$ is a Bernoulli sequence of independent
$\mu$-distributed increments.

For a nested structure $(G,H,X, \phi)$ and a probability measure
$\mu=\mu_G$ on $G$ denote by
$$\mu_X=P_{\phi}(\mu) \hskip .5in \mu_{\bf\til{G}} = \phi(\mu),$$
the images of $\mu$ under the homomorphisms $\phi:G\to \bf\til{G}$
and $P_{\phi}:G \to Sym(X)$.

Then the sample paths of the random walk $({\bf\til{G}},
\mu_{\bf\til{G}})$ (the image of the random walk $(G,\mu)$ under the
embedding $\phi:G\to {\bf\til{G}}$)starting from the identity are
$$({\bf g}_{n}, \si_{n}) = ({\bf h}_{1},\tau_{1})  ({\bf h}_{2},\tau_{2}) \dots  ({\bf h}_{n},\tau_{n}),$$
where $ ({\bf h}_{i},\tau_{i})$ are the
$\mu_{\bf\til{G}}$-distributed increments of the random walk, so
that by definition,
\begin{align}\label{form:rand}
{\bf g}_{n+1} = {\bf g}_{n}{\bf h}^{\si_{n}}_{n+1}, \hskip .5in \si_{n+1}= \si_{n}\tau_{n+1}.\end{align}

\subsection{Reduction to a random walk with internal degrees of freedom.} In particular, formula (\ref{form:rand}) implies that for any fixed $x \in X$
\begin{align} \label{form:ind_rand}
{\bf g}^{x}_{{n+1}} = {\bf g}^{x}_{{n}}\left({\bf h}^{\si_{{n}}}_{{n+1}}\right)^{x} =  {\bf g}^{x}_{{n}}{\bf h}^{x.\si_{{n}}}_{n+1},
\end{align}
so that at any given time $n$ the transition law from ${\bf g}^{x}$
to ${\bf g}_{n+1}^x$ is determined just by the values ${\bf g}^{x}$
and $\si_n(x)$. Therefore the image
\begin{align}\label{form:states}({\bf g}_n^x, x.\si_n)\end{align}
of the original random walk $({\bf \til{G}}, \mu_{\bf\til G})$ under
projection
\begin{align}\label{form:Mark_chain} \Pi_x:{\bf\til{G}} \to H\times X, \hskip .5in ({\bf g}, \si) \to ({\bf g}^x, x.\si)\end{align}
is also a Markov chain.

Formula (\ref{form:ind_rand}) shows that in the notation (\ref{form:notation}) transition of the quotient chains (\ref{form:states}) on
$H\times X$ are
\begin{align}\label{form:trans_prob}(g,z)_{\stackrel{\longrightarrow}{({\bf h}, \tau) \sim\mu_{\bf\til G}}} (g {\bf h}^{z}, z.\tau).\end{align}
Therefore, \\
(i) The transition probabilities (\ref{form:trans_prob}) of the chains (\ref{form:states}) are the same for all $x \in X$;\\
(ii) These transition probabilities are equivariant with respect to
the action of the group $H$ on itself on the left.

Thus the chains (\ref{form:states}) are {\it random walks with internal degrees of
freedom (RWIDF)} on $H$ (the space of these degrees being $X$)
\cite{KS83}.  In other terminologies they are {\it matrix-valued
random walks} \cite{CW89} or {\it covering Markov chains} with the
deck transformation group $H$ and the quotient space $X$ \cite{K95}

Recall that the transition probabilities $$\pi_{g,x}((gh, y)) =
\mu_{x,y}(h)$$ of a general RWIDF are determined by a $d\times d$
matrix
$$M= (\mu_{x,y})_{x,y\in X},$$
(where $d = card X$ is the cardinality of the space of internal
degrees of freedom $X$) of sub-probability measure on the group with
$$\sum_{y} \|\mu_{x,y}\| =1,$$ where $\|\mu_{x,y}\|$ denotes the
total mass of the measure $\mu_{x,y}$. $\mu_{x,y}$ can be treated as
an element of the group algebra.   In the case of a usual random
walk this matrix is $1\times 1$ and consists just of a single
probability measure on the group which determines the random walk.

The image of the RWIDF $(H, M)$ under the map $(g,x) \to x$ is the
quotient Markov chain on $X$ with the transition matrix
\begin{align}\label{form:prob}P=(p_{x,y}) = {\it a}M, \hskip .5in p_{x,y} =\|\mu_{x,y}\|,\end{align}
which is the image of the matrix $M$ under the augmentation
homomorphism ${\it a}: g \to 1$.

In our situation, since $\mu_{\bf \til G}$ is the image of the
measure $\mu$ under the map $\phi$, the matrix $M$ is
\begin{align}\label{form:matrix} M=M^{\mu} = \sum \mu(g) M^g,\end{align}
where $M^g$ are the matrices determined as
$$M_{x,y}^g = \begin{cases} {\bf h}_g^x, \text{ if } & y =x.\si_g \\ 0, &\text{otherwise}\end{cases}.$$

{\bf Remark:} The map $g \to M^g$ is in face an embedding of the
group $G$ into the algebra $Mat(d, \ell^1(H, \rr))$ of $d\times
d$-matrices over the group algebra of $H$.

For the RWIDF $(H, M)$ determined by the matrix $M$ (\ref{form:matrix}) the
quotient chain on $X$ has the transition probabilities (\ref{form:prob})
$$z_{\stackrel{ \longrightarrow}{({\bf h}, \tau) \sim\mu_{\bf \til G}}} z.\tau,$$
or, equivalently
\begin{align}\label{form:chain}z_{\stackrel{\longrightarrow}{\tau \sim \mu_X}} z.\tau.\end{align}

\subsection{Nested reductions.}  Since the Markov chain (\ref{form:chain}) is obtained from
a random walk on the group of permutation $Sym(X)$, any state $x \in
X$ is recurrent. Therefore, the corresponding set $H\times\{x\}$ is
recurrent for the RWIDF (\ref{form:trans_prob}).  Recall that the stopping a Markov
chain at the times when it visits a certain recurrent subset of the
state space gives a new Markov chain on this recurrent subset (the
trace of the original Markov chain).  In our case the transition
probabilities of the induced chain on $H \times \{x\}$ are obviously
equivariant with respect to the left action of the group $H$ on
itself (because the original RWIDF also has this property).
Therefore, the induced chain on $H\times\{x\}$ is actually the usual
random walk determined by a certain probability measure $\mu^x$ on
$H$.

\begin{theorem} \label{thm:induced_meas} The measures $\mu^x$, $x \in X$ can be expressed in terms of the matrix $M$ (\ref{form:matrix}) as
\begin{align}\label{form:ind_meas}\mu^x &= \mu_{x,x} + M_{x, \bar{x}}(I + M_{\bar{x}, \bar{x}}+M_{\bar{x}, \bar{x}}^2 + \dots)
M_{\bar{x}, x} = \\ & =\mu_{x,x} + M_{x, \bar{x}}(I - M_{\bar{x}, \bar{x}})^{-1}  M_{\bar{x}, x},\end{align}
where $M_{x, \bar{x}}$ (resp.,  $M_{\bar{x},x}$ denotes the row
$(\mu_{x,y})_{y\not= x}$ (resp., the column $(\mu_{y,x})_{y\not=x}$)
of the matrix $M$ with the removed element $\mu_{x,x}$ and
$M_{\bar{x}, \bar{x}}$ is the $(d-1)\times(d-1)$ matrix (where $d =
card X$) obtained from $M$ by removing its $x$-th row and column.
The multiplication above is understood in the usual matrix sense in
the group algebra.

\end{theorem}
The proof is the same as in \cite{K04}

\section{Entropy Estimates.}
Recall that the entropy of a probability measure
$\theta=\{\theta_i\}$ is defined as
$$H(\theta) = -\sum \theta_i \log(\theta_i).$$

If $\mu$ is a probability measure on a countable group $G$ with
$H(\mu)<\infty$, then the (asymptotic) entropy of random walk
$(G,\mu)$ is defined as the limit
$$h(G, \mu) = \lim_{n\to \infty} \frac{1}{n}H(\mu_n),$$
where $\mu_n$ denotes the $n$-fold convolution of the measure $\mu$,
i.e., the distribution of the position at time $n$ of the random
walk $(G, \mu)$ issued from the identity of the group $G$.

The following theorem can be easily deduced from Theorem 3.3
\cite{K04}.
\begin{theorem}  Let $(G,H,X,\phi)$ be a nested structure.  Let $\mu$ be a probability measure on $G$ with finite entropy $H(\mu)$.
Then for all $x \in X$
$$h(G, \mu) \leq h(H, \mu^x).$$
\end{theorem}

\subsection{Entropy and amenability.}
Let me quickly recall an amenability criteria.  There are many
definitions of amenability; we will use the Folner condition. Let
$G$ be a countable group.  If for a given finite set $K$ and $\ep>0$
there exists a finite subset $A \ssu G$ such that
$$|A g \Delta A| \leq \ep |A|, \forall g \in K,$$
then we say that $G$ is amenable.

Another equivalent condition of amenability is an existence of a
right-hand invariant mean on $G$.

The fundamental fact relating the asymptotic entropy with
amenability of the group $G$ is that $h(G, \mu)=0$ if and only if
the Poisson boundary of the random walk $(G, \mu)$ is trivial. This
in turn implies amenability (as the action of the group on Poisson
boundary is amenable). Moreover $n$-fold convolution of the measure
$\mu$ converges to invariant mean \cite{KV83}.  Therefore, if a
group $G$ carries a non-degenerate random walk with vanishing
asymptotic entropy, then it must be amenable.

\section{Grigorchuk $2$-groups.}
In this section I describe a construction of Grigorchuk $2$-group.
For complete description and further results see \cite{Gr84}

Let $\Om = \{0,1,2\}^{\nn}$ and $\si: \Om\to \Om$ is left side
shift.  For every $\om \in \Om$, there exists the Grigorchuk
$2$-group $G_{\om}$ which is a subgroup of the group of
automorphisms of the binary tree $Aut(T_2)$.  The group $G_{\om}$ is
generated by $4$ elements $a, b_{\om}, c_{\om}, d_{\om} \in
Aut(T_2)$.

To describe the group I need to specify the action of each element
on the binary tree.  
Each vertex of the tree can be represented as a word in the $2$-letter alphabet $\{0,1\}$.
The action on
the tree  is easy to recover from that definition.

Define $t_i:\{0,1,2\} \to Aut(T_2)$ as
$$t_i(j) = \begin{cases} 1 &\text{ if } i=j \\ a &\text{ if } i\not=j \end{cases}.$$
For $x \in \{0,1\}^{\nn}$ I have

\begin{align*} &a(0x) = 1x, &a(1x) = 0x \\
& b_{\om}(0x)=0t_0(\om_1)(x), &b_{\om}(1x)= 1b_{\si(\om)}(x),\\
 & c_{\om}(0x)=0t_1(\om_1)(x), &c_{\om}(1x)= 1c_{\si(\om)}(x), \\
&d_{\om}(0x)=0t_2(\om_1)(x), &d_{\om}(1x)=
1d_{\si(\om)}(x),\end{align*} where $\om =
(\om_1,\om_2,\om_3,\dots)$.

Let $\varepsilon:\{0,1\}\to \{0,1\}$ is defined as $\varepsilon(i) =
1-i$. Observe that there exists $\pi_i: G_{\om} \to G_{\si \om}$ for
$i\in \{0,1\}$, where $\pi_i$ are defined as $g(ix) =
\varepsilon^{l(g)}(i)\pi_i(g)(x),$ where $ l(g) \in \{0,1\}$. This
proves
\begin{lemma} \label{lem:nested_str} The embedding $$\phi_{\om}: G_{\om} \to G_{\si \om}^{\{0,1\}} \leftthreetimes Sym(\{0,1\}),$$
defined as $g \to (\pi_0(g), \pi_1(g), \varepsilon^{l(g)})$, induces
the nested structure
$$(G_{\om}, G_{\si \om}, \{0,1\}, \phi_{\om}).$$
\end{lemma}

Let me give another description of $G_{\om}$ that is easier is to
visualize.  It is the original Grigorchuk description of $G_{\om}$.
The space of ends is ''almost'' the interval $[0,1]$.

Let $\Delta$ be an interval.  Denote by $I$ an identity
transformation on
 $\Delta$
and by $T$ a transposition of two halves of $\Delta$.

For each $\omega \in \Omega$ define a $3 \times \infty$ matrix
${\overline{\omega}}$ by replacing $\omega_{i}$ with columns
$\overline{\omega}_i$ where
$$\bar {0}=
\left(
\begin{array}{c}
T\\
T\\
I\\
\end{array}
\right)\ ,\ \bar {1}= \left(
\begin{array}{c}
T\\
I\\
T\\
\end{array}
\right)\ ,\ \bar {2}= \left(
\begin{array}{c}
I\\
T\\
T\\
\end{array}
\right)
$$

By $U^{\omega}= (u_1^{\omega},u_2^{\omega},\dots),
V^{\omega}=(v_1^{\omega},v_2^{\omega},\dots),
W^{\omega}=(w_1^{\omega},w_2^{\omega},\dots)$ denote the rows of
${\overline{\omega}}$. Think of them as of infinite words in the
alphabet $\{T,I\}$.

Define transformations $a_{\omega},b_{\omega},c_{\omega},d_{\omega}$
of an interval $\De = [0,1]\sm  \mathbb{Q}$ as follows:

\medskip

\begin{picture}(400,100)(0,0)
\put (10,80){$a_\omega :$} \put (40,80){\line(1,0){100}} \put
(40,70){$0$} \put (140,70){$1$} \put (90,85){$T$}

\put (10,40){$b_\omega :$} \put (40,40){\line(1,0){100}} \put
(40,30){$0$} \put (140,30){$1$} \put (87,37){$\circ$} \put
(87,27){$1\over 2$} \put (112,37){$\circ$} \put (112,27){$3\over4$}
\put (120,30){$\dots$} \put (65,45){$u_1^{\omega}$} \put
(100,45){$u_2^{\omega}$} \put (120,45){$\dots$}

\put (210,80){$c_\omega :$} \put (240,80){\line(1,0){100}} \put
(240,70){$0$} \put (340,70){$1$} \put (287,77){$\circ$} \put
(287,67){$1\over 2$} \put (312,77){$\circ$} \put (312,67){$3\over4$}
\put (320,70){$\dots$} \put (265,85){$v_1^{\omega}$} \put
(300,85){$v_2^{\omega}$} \put (320,85){$\dots$}

\put (210,40){$d_\omega :$} \put (240,40){\line(1,0){100}} \put
(240,30){$0$} \put (340,30){$1$} \put (287,37){$\circ$} \put
(287,27){$1\over 2$} \put (312,37){$\circ$} \put (312,27){$3\over4$}
\put (320,30){$\dots$} \put (265,45){$w_1^{\omega}$} \put
(300,45){$w_2^{\omega}$} \put (320,45){$\dots$}
\end{picture}

\medskip

Observe that $a_{\omega}$ is independent of $\omega$, and will be
further denoted by $a$. Let $G_{\omega}$ be a group of
transformations of the interval $\Delta$ generated by $a,
b_{\omega}, c_{\omega}, d_{\omega}$.

It was proved by Grigorchuk \cite{Gr84} that if $\om$ does not
become constant, (i.e., there does not exists $N$ such that
$\om_{N+1}=\om_{N+2}= \dots$) the group $G_{\om}$ has subexponential
growth. In case when there exists $N$ such that
$\om_N\not=\om_{N+1}=\om_{N+2}= \dots$, $G_{\om}$ has a subgroup of
finite index isomorphic to the free abelian group of rank $2^N$

Before I proceed to the next section I mention the result prove by
Grigorchuk \cite{Gr84} (see Theorem 7.1)

\begin{theorem}
Suppose that function $\rho(n)$ grows more slowly that any
exponential function, i.e., $\rho(n)=o(2^{\ep n})$ for all $\ep>0$.
Then there exist $\om \in \Om$ and $N$ such that
$|B_{\om}(N)|>\rho(N),$ where $B_{\om}(n)$ is a ball of radius $n$
in $G_{\om}$ with respect to generators $a, b_{\om}, c_{\om},
d_{\om}$.
\end{theorem}

{\bf Remark:} This is not the exact theorem 7.1 from \cite{Gr84}.
But it is a simplified version of it.

\section{Constructing the universal group}

Let $X$ be a set and $(G, \odot)$ be a group.  Denote by $(G^X,
\star)$ be the group of all maps from $X$ to $G$ with pointwise
composition.

\begin{proposition} Let  $\si$ be a surjective map from $X$ to $X$.
Then $$\hat{\si}: G^X \to G^X$$ defined as $f \to f\circ \si$ is an
injective homomorphism.
\end{proposition}

Before I proceed I will need a simple lemma.

\begin{lemma} Let $f : X\to X$ be a surjective map.  Assume $g:X \to Y$ be some
map.  Then $g\equiv const$ if and only if $g \circ f \equiv const$
\end{lemma}

$\square$ Since $f$ is surjective, there exists right inverse, i.e.,
$h:X \to X$ such that $f \circ h = Id$.  Therefore if $g \circ f
\cong const$ then $g = (g \circ f) \circ h \equiv const \circ h =
const$. $\blacksquare$

{\bf Proof of Proposition:}  Let $f, g \in G^X$, then
$(\hat{\si}(f\star g))(x) = (f\star g)(\si(x)) = f(\si(x))\odot
g(\si(x)) = \hat{\si}(f)(x) \odot \hat{\si}(g)(x) = (\hat{\si}(f)
\star \hat{\si}(g))(x)$ for all $x \in X$. Since $Id\in G^X$ is a
constant map, $\hat{\si}(Id)=Id$. Thus $\hat{\si}$ is a
homomorphism.

Now if $f \in ker(\hat{\si})$ then $ f \circ \si = Id$. By the above
lemma $f =Id$.  Thus $\hat{\si}$ is injective. $\blacksquare$

Fix $\phi_1, \phi_2, \dots, \phi_n \in G^X$.  Then for each $\omega
\in X$, I define $$H_{\omega} = \< \phi_1(\omega), \dots,
\phi_n(\omega)\> \leq G.$$

It is easy to observe that $H = \< \phi_1, \phi_2, \dots, \phi_n\>$
is the smallest group such that for every $\omega \in X$,
$H_{\omega}$ is a quotient of $H$.

\section{Universal Grigorchuk $2$-group.}
Let $X = \Om=\{0,1,2\}^{\zz_+}$, $H= Aut(T_2)$, where $T_2$ is a
binary tree, and $\si$ is a shift on $\Om$.

I am going to consider a subgroup of $H^{\Om}$ generated by maps
$A = a$, where $a$ is a flip at the top vertex. $B(\om) = b_{\om}$,
$C(\om) = c_{\om}$, $D(\om) = d_{\om}$, where
$b_{\om},c_{\om},d_{\om}$ are element of the Grigorchuk group
corresponding to $\om$.

Let $Gr_2 = \< A, B, C, D\>$.  I will call $Gr_2$ the universal
Grigorchuk $2$-group.  It is clear that for every $\om \in \Om$, the
evaluation map at $\om$ is a surjective map form $G$ to $G_{\om}$.

\begin{lemma} The group $Gr_2$ has exponential growth. \end{lemma}

$\square$  Let $\rho(n)$ be the growth function with respect to
generators $A, B, C,D$.  Assume toward a contradiction, that $Gr_2$
does not have exponential growth. Then $\lim_{n \to \infty}
\frac{\log(\rho(n))}{n} =0$.  Thus $\rho(n)$ has subexponential
growth, i.e., $\rho(n) \leq e^{\epsilon n}$ for every $\epsilon>0$.
By Theorem 7.1 in \cite{Gr84}, there exists $\om$ and $N$ such that
$B_{\om}(n) > \rho(n)$, where $B_{\om}(n)$ is the size of the ball
of radium $n$ inside the group $G_{\om}$ with respect to generators
$a_{\om}, b_{\om}, c_{\om}, d_{\om}$.

However $G_{\om}$ is a quotient of $G$ and generators $A, B, C,D$
maps to $a_{\om}, b_{\om}, c_{\om}, d_{\om}$. Thus for every $n$,
$\rho(n) \geq B_{\om}(n)$.  But for $N$, I have $\rho(N) \geq
B_{\om}(N)>\rho(N).$ Contradiction.  $\blacksquare$

\begin{lemma} The map $\phi: Gr_2 \to \si(Gr_2)^{\{0,1\}} \leftthreetimes Sym(\{0,1\})$ defined
by $$\phi(A) = (Id, Id, \varepsilon),$$  $$\phi(B) = (\til{t}_0,
B\circ \si, Id),$$
$$\phi(C) = (\til{t}_1, C\circ \si, Id), $$
$$\phi(D) = (\til{t}_2, D\circ \si, Id)$$
induces a nested structure $(Gr_2, \si(Gr_2), \{0,1\}, \phi)$, where
$$\til{t}_i( \om) = \begin{cases} Id & \text{ if } \om_1=i \\a  &
\text{ if } \om_1\not=i \end{cases}.$$
\end{lemma}

$\square$ The proof is an easy consequence of Lemma
\ref{lem:nested_str}. The hard part is to show that is it an
embedding.  If $W(A,B,C,D)$ is a word such that $\phi(W(A,B,C,D))=
(Id,Id, Id)$, then $\phi_{\om}(W(A,B,C,D)(\om))= (Id,Id, Id)$ for
all $\om \in \Om$.  But since $\phi_{\om}$ is an embedding I have
$$W(A(\om),B(\om),C(\om),D(\om)) = Id$$ for every $\om \in \Om$.
This proves the lemma. $\blacksquare$
\section{Amenability of Universal Grigorchuk $2$-group.}

In this section I complete the proof of Theorem \ref{thm:main}.

$\square$ Let $Gr_2 = \<A, B, C, D\>$ then $Gr_2 \cong
\hat{\si}(Gr_2) = \<A, B\circ \si, C\circ \si, D\circ \si\>$, as $A
\circ \si = A$.

Since we have a nested structure, we can write everything in matrix
form.

$A \leftrightarrow \begin{pmatrix} 0 & Id \\ Id & 0\end{pmatrix}$,
$B \leftrightarrow \begin{pmatrix} \til{t}_0 & 0 \\ 0 & B \circ \si
\end{pmatrix}$, $C \leftrightarrow \begin{pmatrix} \til{t}_1 & 0 \\
0 & C \circ \si \end{pmatrix}$, $D \leftrightarrow \begin{pmatrix}
\til{t}_2 & 0 \\ 0 & D \circ \si \end{pmatrix}.$

A trivial observation is that $\til{t}_1 + \til{t}_2 + \til{t}_3 =
2A + 1 \in R(G) \ssu R(H^X)$ in the group algebras. Consider a
random walk on $Gr_2$ generated by $\mu_0=\alpha Id + \beta A + m
(B+C+D)$. This corresponds to

$$ \mu_0 \leftrightarrow \begin{pmatrix}
\alpha Id + m (2A + Id) & \beta Id \\ \beta Id & \alpha Id+ m(B\circ \si + C \circ \si + D\circ \si)\end{pmatrix}.$$

Consider an induced random walk on the left subtree.  (I will denote
induced measure with superscript $ind$.) 
It will be supported on $A,
B\circ \si , C \circ \si , D\circ \si$ and using formula (\ref{form:ind_meas}) from Theorem \ref{thm:induced_meas} we have
\begin{align*}
\mu_1&=\mu_0^{ind}=(\alpha + m) Id + 2m A + \beta^2 \sum_{i=0}^{\infty} \left(\alpha Id+ m(B\circ \si + C \circ \si + D\circ
\si)\right)^i= \\ &=
(\alpha + m + x) Id + 2m A + \frac{\beta -
x}{3}(B\circ \si +C\circ \si+D\circ \si),\end{align*} where $x$ is some
non-negative number depending on $\beta$ and $m$, coming from the
term
$$\beta^2 \sum_{i=0}^{\infty} \left(\alpha Id+ m(B\circ \si + C \circ \si + D\circ
\si)\right)^i,$$
and using the fact that $B,C,D$ commute and $B^2=C^2=D^2=BCD=Id$.

{\bf Remark:} I have used that $\mu_0$ and $\mu_1$ are probabilities measures to make the calculations.

Now recall that $\hat{\si}(Gr_2) \cong Gr_2$. Therefore, this gives
a random walk on the initial group.  It is possible to calculate
$\mu_1$ explicitely and then find $\mu_0$ that satisfies the
equation $\mu_0=\delta Id + (1-\delta)\mu_1$ for some $0<\delta<1$.
In this case the problem would be reduced to the method of V.
Kaimanovich in \cite{K04}.

I would like to use a small trick in order to avoid calculations of
$x$. We can induce a random walk on the left subtree, again giving
rise to a new measure

$$\mu_2=\mu_1^{ind}=(\alpha + m + x + \frac{\beta - x}{3} + y) Id +
 2\frac{\beta - x}{3} A + \frac{2m - y}{3}(B\circ \si^2 +C\circ \si^2+D\circ \si^2),$$

Now it is easy to observe that the coefficients near $A, B, C, D$
have shrunk by a factor at least $2/3$. Since $\mu_n$ are probability measures,  $\mu_n \to
Id$
as $n \to \infty$, where
$\mu_{n+1}=\mu_n^{ind}$. As $h(\mu) \leq H(\mu)$, and $h(\mu_i) \leq
h(\mu_{i+1})$ it follows that
$$h(\mu_0) \leq limsup_{n \to \infty} h(\mu_n) \leq lim_{n\to \infty} H(\mu_n) =0.$$

This proves that $Gr_2$ is amenable. $\blacksquare$

{\bf Remark:} In the last part of the proof, it is possible to
actually calculate $x$ and then construct a self-similar measure as
in \cite{K04}.  I did not do that in order to show another way to
estimate the entropy of the random walk.

\bibliographystyle{alpha}

\end{document}